\documentclass[11pt, a4paper]{article}
\usepackage{latexsym}
\usepackage{indentfirst}
\usepackage{graphicx}
\usepackage{amsmath}
\usepackage{amssymb}

\begin{document}
\title{Approximation of pressure perturbations by FEM}
\author{C\u{a}t\u{a}lin Liviu Bichir $^{1}$, Adelina
Georgescu $^{2}$, \\
1 \ Research and Design Institute for Shipbuilding - ICEPRONAV \\
Gala\c{t}i, Romania, \\
2 \ University of Pite\c{s}ti, Romania}

\maketitle

\begin{abstract}
In the mathematical problem of linear hydrodynamic stability for
shear flows against Tollmien-Schlichting perturbations, the
continuity equation for the perturbation of the velocity is
replaced by a Poisson equation for the pressure perturbation. The
resulting eigenvalue problem, an alternative form for the two -
point eigenvalue problem for the Orr - Sommerfeld equation, is
formulated in a variational form and this one is approximated by
finite element method (FEM).
Possible applications to concrete cases are revealed. 
\end{abstract}

\section{The classical model and the model based on
a Poisson equation for the pressure perturbation}
\label{sectiunea1}

The mathematical problem of linear hydrodynamic stability is
\cite{AG200}
\begin{subequations}
\label{NS1}
\begin{eqnarray}
   & \ & {\partial \textbf{u}' \over \partial t}
            + (\textbf{U}\cdot \textbf{grad})\textbf{u}'
            + (\textbf{u}'\cdot \textbf{grad})\textbf{U} =
            \label{NS1a} \\
   & \ & \qquad\qquad  = -\textbf{grad} \ p'
            + Re^{-1}\triangle \textbf{u}', \ (t,\textbf{x})\in \mathbb{R}_{+}\times\Omega,
            \nonumber  \\
   & \ & div \ \textbf{u}' = 0, \ (t,\textbf{x})\in \mathbb{R}_{+}\times\Omega,
            \label{NS1b} \\
   & \ & \textbf{u}'(t,\textbf{x}) = 0, \ t\in \mathbb{R}_{+}, \ \textbf{x}\in\partial\Omega,
            \label{NS1c} \\
   & \ & \textbf{u}'(0,\textbf{x}) = \textbf{u}'_{0}(\textbf{x}), \ \textbf{x}\in\Omega,
            \label{NS1d}
\end{eqnarray}
\end{subequations}
and its solution is the perturbation $(\textbf{u}',p')$, for
$t>0$. Here all variables are non dimensional, $(\textbf{U},P)$ is
the basic motion in the domain $\Omega \subseteq \mathbb{R}^{n}$,
$n = 2$ or $3$, $Re=U_{\infty}L / \nu$ is the Reynolds number, $L$
is a length scale, $U_{\infty}$ is a velocity scale, $\nu$ is the
coefficient of kinematic viscosity.

Applying the divergence operator to (\ref{NS1a}) and using the
divergence-free condition on $\textbf{u}'$, we obtain
\begin{equation}
\label{NS1_P_EC}
\begin{split}
   &   \triangle p' =
          - div \cdot ((\textbf{U}\cdot \textbf{grad})\textbf{u}')
          - div \cdot ((\textbf{u}' \cdot \textbf{grad})\textbf{U}), \   \\
   &   \qquad\qquad  (t,\textbf{x})\in \mathbb{R}_{+}\times\Omega,   \\
\end{split}
\end{equation}
and projecting (\ref{NS1a}) on $\textbf{n}_{\ast}$, we have
\begin{equation}
\label{NS1_P_CL}
\begin{split}
   &   \textbf{n}_{\ast}(\textbf{x}) \cdot \textbf{grad} \, p'(t,\textbf{x})
          = \textbf{n}_{\ast}(\textbf{x}) \cdot
          \textbf{F}(Re,\textbf{U}(\textbf{x}),\textbf{u}'(t,\textbf{x})), \   \\
   &   \qquad\qquad
          \ t\in \mathbb{R}_{+}, \ \textbf{x}\in\partial\Omega,
\end{split}
\end{equation}
where $\textbf{n}_{\ast}$ is the unit outward normal to
$\partial\Omega$ and $\textbf{F}(Re,\textbf{U},\textbf{u}')$ $=$
$Re^{-1} \triangle \textbf{u}'$ $-$ $(\textbf{U}\cdot
\textbf{grad})\textbf{u}'$ $-$ $(\textbf{u}' \cdot
\textbf{grad})\textbf{U}$. Another form of the right hand - side
of the Poisson equation for $p'$ can be obtained by using the
divergence-free condition on $\textbf{U}$.

Thus the mathematical problem of the linear hydrodynamic stability
becomes (\ref{NS1a}), (\ref{NS1_P_EC}), (\ref{NS1c}),
(\ref{NS1_P_CL}), (\ref{NS1d}). Of course, it is assumed that
equation (\ref{NS1a}) can be continued on the $\partial\Omega$.

\section{An alternative approximation of the two - point
eigenvalue problem for Orr - Sommerfeld equation}
\label{sectiunea2}

Let $\Omega=\{\textbf{x} = (x,y,z)\in\mathbb{R}^{3} | -\infty
<x,z< \infty, 0<y<a \}$ ($a \geq 1$) and assume that the basic
flow is of the form $\textbf{U}(x,y,z)=(U(y),0,0)$. Let us choose
Tollmien-Schlichting waves - like perturbations $\textbf{u}'(t,x,
y, z)=\textbf{u}'_{0}(y)exp[i \alpha (x-ct)]$,
$\textbf{u}'_{0}(y)=(u(y),v(y),0)$, $\textbf{u}'(t,x, y, z)$ $=$
$(u'(t,x, y),v'(t,x, y),0)$, $p'(t,x, y, z) =p(y)exp[i \alpha
(x-ct)]$, where $i=\sqrt{-1}$, $\alpha$ is the streamwise wave
number, $c=c_{r} + i c_{i}$, $c_{r}$ is the wave speed and $c_{i}$
is the amplification rate.

In this case, model (\ref{NS1}) leads to the classical two - point
eigenvalue problem for Orr - Sommerfeld equation in $(c,\varphi)$,
where $\varphi$ is the nonexponential factor of the stream
function, while (\ref{NS1a}), (\ref{NS1_P_EC}), (\ref{NS1c}),
(\ref{NS1_P_CL}), (\ref{NS1d}) becomes the two - point eigenvalue
problem in $(c,u,v,p)$
\begin{subequations}
\label{NS1_C}
\begin{eqnarray}
   & \ & - \alpha Re^{-1} i u
          + (\alpha Re)^{-1} i u''
          + U u
          - \alpha^{-1} i U' v
          + p
          = c u, \ y \in (0,a), \label{NS1_Ca} \\
   & \ & - \alpha Re^{-1} i v
          + (\alpha Re)^{-1} i v''
          + U v
          - \alpha^{-1} i p'
          = c v, \ y \in (0,a), \label{NS1_Cb} \\
   & \ & - \alpha^{2} p
          + p''
          = - 2 \alpha i U' v, \ y \in (0,a), \label{NS1_Cc} \\
   & \ & u(y) = 0
          \ \ \textrm{for} \ y=0, y=a, \label{NS1_Cd} \\
   & \ & v(y) = 0
          \ \ \textrm{for} \ y=0, y=a, \label{NS1_Ce} \\
   & \ & p'(y)= Re^{-1} v''(y)
          \ \ \textrm{for} \ y=0, y=a. \label{NS1_Cf}
\end{eqnarray}
\end{subequations}
where the prime stands for the differentiation with respect to
$y$, $c$ is an eigenvalue and $(u,v,p)$ is an eigenvector.

Let $L^{2}(0,a)$ be the Hilbert space of all measurable
complex-valued functions u, defined on $(0,a)$, for which
$\int\limits_{0}^{a}|u(x)|^{2}dx<\infty$ and the inner product is
$(u,v)=\int\limits_{0}^{a}u\overline{v} dy$. Denote $Du$ the
generalized derivative and consider the spaces
$H_{0}^{1}(0,a)=\{u\in L^{2}(0,a)|Du \in L^{2}(0,a), u(0)=u(a)=0
\}$, $L_{0}^{2}(0,a) = \{ p \in L^{2}(0,a) \ | \
\int\limits_{0}^{a}p(x) \ dx = 0 \}$.

Multiply (\ref{NS1_Ca}), (\ref{NS1_Cb}) and (\ref{NS1_Cc}) by
arbitrary test functions $\overline{f}_{1}$, $\overline{f}_{2}$
and $\overline{g}$ respectively, integrate over $(0,a)$, apply
partial integration if necessary and take into account the two
point conditions (\ref{NS1_Cd}) - (\ref{NS1_Cf}) to obtain the
following weak formulation of problem (\ref{NS1_C}) in $(c,u,v,p)$
$\in$ $\mathbb{C}$ $\times$ ${H_{0}^{1}(0,a)}^{2}$ $\times$
$L_{0}^{2}(0,a)$
\begin{equation}
\label{eq_0}
\begin{split}
   & - \alpha Re^{-1} i \int\limits_{0}^{a}u \overline{f}_{1} dy
        - (\alpha Re)^{-1} i \int\limits_{0}^{a}
           u' {\overline{f}_{1}}' dy
        + \int\limits_{0}^{a}U u \overline{f}_{1} dy -  \\
   & \quad
        - \alpha^{-1} i \int\limits_{0}^{a}U' v \overline{f}_{1} dy
        + \int\limits_{0}^{a} p \overline{f}_{1} dy
        = c \cdot \int\limits_{0}^{a} u \overline{f}_{1} dy,
        \forall \ f_{1} \in H_{0}^{1}(0,a),
\end{split}
\end{equation}

\begin{equation}
\label{eq_1}
\begin{split}
   & - \alpha Re^{-1} i \int\limits_{0}^{a}v \overline{f}_{2} dy
        - (\alpha Re)^{-1} i \int\limits_{0}^{a}
           v' {\overline{f}_{2}}' dy
        + \int\limits_{0}^{a}U v \overline{f}_{2} dy  \\
   & \quad
        + \alpha^{-1} i \int\limits_{0}^{a} p {\overline{f}_{2}}' dy
        = c \cdot \int\limits_{0}^{a} v \overline{f}_{2} dy,
        \forall \ f_{2} \in H_{0}^{1}(0,a),
\end{split}
\end{equation}

\begin{equation}
\label{eq_2}
\begin{split}
   & - \alpha^{2} \int\limits_{0}^{a} p \overline{g} dy
        - \int\limits_{0}^{a}
           p' {\overline{g}} \: ' dy
        = - Re^{-1} v''(a) \overline{g}(a)
           \\
   & \quad
        + Re^{-1} v''(0) \overline{g}(0)
        - 2 \alpha i \int\limits_{0}^{a}U' v \overline{g} dy,
        \forall \ g \in L_{0}^{2}(0,a).
\end{split}
\end{equation}

\section{Approximation of problem (\ref{eq_0}) -
(\ref{eq_2}) by FEM} \label{sectiunea3}

In order to perform this approximation, let us divide the interval
$[0,a]$ in $N+1$ subintervals $K=K_{j}=[y_{j},y_{j+1}]$, $0 \leq j
\leq N$, where $0 = y_{0} < y_{1} < \ldots < y_{N+1} = a$. The
sets $K$ represent a triangulation $\mathcal{T}_{h}$ of $[0,a]$.

The approximate basic shear flow
$\textbf{U}_{h}(x,y,z)=(U_{h}(y),0,0)$ is defined by
$U_{h}(y)=U(y)$, $y \in [0,a]$. The approximate amplitudes
$u_{h}(y)$, $v_{h}(y)$ and $p_{h}(y)$ correspond to the exact
ones, $u(y)$, $v(y)$ and $p(y)$, respectively.

Introduce the spaces $V_{h} = \{ v : [0,a] \rightarrow \mathbb{C}
\: | \: v \in C[0,a]$, $v$ is on $K_{j}$ an one - dimensional
polynomial, having complex coefficients, of degree $2$, $0 \leq j
\leq N, v(0) = v(a) = 0 \}$ and $M_{h} = \{ v : [0,a] \rightarrow
\mathbb{C} \: | \: v \in C[0,a]$, $v$ is on $K_{j}$ an one -
dimensional polynomial, having complex coefficients, of degree
$1$, $0 \leq j \leq N \}$.

Correspondingly, variational problem (\ref{eq_0}) - (\ref{eq_2})
is approximated by the following problem in
$(c_{h},u_{h},v_{h},p_{h})$ $\in$ $\mathbb{C}$ $\times$
$V_{h}^{2}$ $\times$ $M_{h}$,

\begin{equation}
\label{eq_0h}
\begin{split}
   & - \alpha Re^{-1} i \int\limits_{0}^{a}u_{h} \overline{f}_{1} dy
        - (\alpha Re)^{-1} i \int\limits_{0}^{a}
           u'_{h} {\overline{f}_{1}}' dy
        + \int\limits_{0}^{a}U u_{h} \overline{f}_{1} dy -  \\
   & \quad
        - \alpha^{-1} i \int\limits_{0}^{a}U' v_{h} \overline{f}_{1} dy
        + \int\limits_{0}^{a} p_{h} \overline{f}_{1} dy
        = c_{h} \cdot \int\limits_{0}^{a} u_{h} \overline{f}_{1} dy,
        \forall \ f_{1} \in V_{h},
\end{split}
\end{equation}

\begin{equation}
\label{eq_1h}
\begin{split}
   & - \alpha Re^{-1} i \int\limits_{0}^{a}v_{h} \overline{f}_{2} dy
        - (\alpha Re)^{-1} i \int\limits_{0}^{a}
           v'_{h} {\overline{f}_{2}}' dy
        + \int\limits_{0}^{a}U v_{h} \overline{f}_{2} dy  \\
   & \quad
        + \alpha^{-1} i \int\limits_{0}^{a} p_{h} {\overline{f}_{2}}' dy
        = c_{h} \cdot \int\limits_{0}^{a} v_{h} \overline{f}_{2} dy,
        \forall \ f_{2} \in V_{h},
\end{split}
\end{equation}

\begin{equation}
\label{eq_2h}
\begin{split}
   & - \alpha^{2} \int\limits_{0}^{a} p_{h} \overline{g} dy
        - \int\limits_{0}^{a}
           p'_{h} {\overline{g}} \: ' dy
        = - Re^{-1} v''_{h}(a) \overline{g}(a)
           \\
   & \quad
        + Re^{-1} v''_{h}(0) \overline{g}(0)
        - 2 \alpha i \int\limits_{0}^{a}U' v_{h} \overline{g} dy,
        \forall \ g \in M_{h}.
\end{split}
\end{equation}

In order to obtain $u_{h}$, $v_{h}$, we use a basis of real
functions of $V_{h}$. Let $J=J_{K}=\{ 1,2,3 \}$ be the local
numeration for the nodes of $K$, where $1, 3$ correspond to
$y_{j}, y_{j+1}$ respectively and $2$ corresponds to a node
between $y_{j}$ and $y_{j+1}$. Let $\{ \phi_{n}, n \in J \}$ be
the local quadratic basis of functions on $K$ corresponding to the
local nodes. Let $J_{\ast}=\{ 1, \ldots, 2N+1 \}$ be the global
numeration for the nodes of $[0,a]$, where the nodes corresponding
to $y_{0}$ and $y_{N+1}$ are not taken into account, and let
$L_{1}$ be a matrix whose elements are the elements of $J_{\ast}$.
Its rows are indexed by the elements $K \in \ \mathcal{T}_{h}$ and
its columns, by the local numeration $n \in \ J$. We take the
value $0$ at the locations of $L_{1}$ which we do not consider in
the computations, i.e. the locations where $K=K_{0}$, $n=1$ and
$K=K_{N}$, $n=3$. Write $n_{\ast} = L_{1}(K,n)$ or, simply,
$n_{\ast}$ for the element $n_{\ast}$ of $L_{1}$ which depends on
$K$ and $n$. Let $\{ (\Phi_{n_{\ast}},0),(0,\Phi_{n_{\ast}}) ;
n_{\ast} \in J_{\ast} \}$ be a basis of functions of $V_{h}^{2}$.
If $n_{\ast}$ corresponds to $y_{j}$, then $\Phi_{n_{\ast}}$ is a
real quadratic function on $K_{j-1}$ and $K_{j}$ and its value is
zero on $[0,a] \backslash (K_{j-1} \cup K_{j})$. If $n_{\ast}$
lies between $y_{j}$ and $y_{j+1}$, then $\Phi_{n_{\ast}}$ is a
real quadratic function on $K_{j}$ and its value is zero on $[0,a]
\backslash K_{j}$. We have $\Phi_{n_{\ast}}(y) = \phi_{n}(y)$,
where $y \in K$, $n_{\ast} = L_{1}(K,n)$. Let $u_{n_{\ast}}$,
$v_{n_{\ast}}$ be the values of $u_{h}$, $v_{h}$ at the nodes
$n_{\ast}$, $n_{\ast} \in J_{\ast}$. Retaining our convention
about $n_{\ast} = L_{1}(K,n)$, we do not write in the sequel the
conditions $n \neq 1$ for $K=K_{0}$ and $n \neq 3$ for $K=K_{N}$.
We have
\begin{equation*}
\label{eq_10}
   u_{h}(y)
      = \sum\limits_{n_{\ast} \in J_{\ast}}
         u_{n_{\ast}}\Phi_{n_{\ast}}(y)
      = \sum_{K \in \ \mathcal{T}_{h}, y \in K}
        \sum_{n \in \ J}
        u_{n_{\ast}}\phi_{n}(y)
\end{equation*}
and a similar expression for $v_{h}(y)$.

In order to obtain $p_{h}$, we use a basis of real functions of
$M_{h}$. Let $I=I_{K}=\{ 1,2 \}$ be the local numeration for the
nodes of $K$, where $1, 2$ correspond to $y_{j}, y_{j+1}$
respectively. Let $\{ \psi_{m}, m \in I \}$ be the local affine
basis of functions on $K$ corresponding to the local nodes. Let
$I_{\ast}=\{ 0, 1, \ldots, N, N+1 \}$ be the global numeration for
the nodes of $[0,a]$ and let $L_{2}$ be a matrix whose elements
are the elements of $I_{\ast}$. Its rows are indexed by the
elements $K \in \ \mathcal{T}_{h}$ and its columns, by the local
numeration $m \in \ I$. Write $m_{\ast} = L_{2}(K,m)$ or, simply,
$m_{\ast}$ for the element $m_{\ast}$ of $L_{2}$ which depends on
$K$ and $m$. Let $\{ \Psi_{m_{\ast}}, m_{\ast} \in I_{\ast} \}$ be
a basis of functions of $M_{h}$. If $m_{\ast}$ corresponds to
$y_{j}$, then $\Psi_{m_{\ast}}$ is a real affine function on
$K_{j-1}$ and $K_{j}$ and its value is zero on $[0,a] \backslash
(K_{j-1} \cup K_{j})$. We have $\Psi_{m_{\ast}}(y) = \Psi_{m}(y)$,
where $y \in K$, $m_{\ast} = L_{2}(K,m)$. Let $p_{m_{\ast}}$ be
the values of $p_{h}$ at the nodes $m_{\ast}$, $m_{\ast} \in
I_{\ast}$. We have
\begin{equation*}
\label{eq_11}
   p_{h}(y)
      = \sum\limits_{m_{\ast} \in I_{\ast}}
        p_{m_{\ast}}\Psi_{m_{\ast}}(y)
      = \sum_{K \in \ \mathcal{T}_{h}, y \in K}
        \sum_{m \in \ I}
        p_{m_{\ast}}\psi_{m}(y).
\end{equation*}

Problem (\ref{eq_0h}) -  (\ref{eq_2h}) becomes the following
algebraic eigenvalue problem in $(c_{h},$ $u_{1},v_{1},$ $\ldots,$
$u_{2N+1},v_{2N+1},$ $p_{0},$ $\ldots,$ $p_{N+1})$ $\in$
$\mathbb{C}^{1+2(2N+1)+N+2}$
\begin{equation}
\label{eq_20}
\begin{split}
   & \sum_{K \in \ \mathcal{T}_{h}}
        \sum_{n \in \ J} u_{n_{\ast}}[
           - \alpha Re^{-1} i \int\limits_{K}\phi_{n}\phi_{k}dz
           - (\alpha Re)^{-1} i \int\limits_{K}\phi'_{n}\phi'_{k}dz +  \\
   & \quad
           + \int\limits_{K}U\phi_{n}\phi_{k}dz]
           + \sum_{K \in \ \mathcal{T}_{h}}
              \sum_{n \in \ J} v_{n_{\ast}}[
                 - \alpha^{-1} i \int\limits_{K}U'\phi_{n}\phi_{k}dz] +  \\
   & \quad
        + \sum_{K \in \ \mathcal{T}_{h}} \sum_{m \in \ I} p_{m_{\ast}}
           \int\limits_{K}\psi_{m}\phi_{k}dz
        = c_{h} \cdot \{ \sum_{K \in \ \mathcal{T}_{h}} \sum_{n \in \ J} u_{n_{\ast}}
           \int\limits_{K}\phi_{n}\phi_{k}dz \}
\end{split}
\end{equation}

\begin{equation}
\label{eq_21}
\begin{split}
   & \sum_{K \in \ \mathcal{T}_{h}}
        \sum_{n \in \ J} v_{n_{\ast}}[
           - \alpha Re^{-1} i \int\limits_{K}\phi_{n}\phi_{k}dz
           - (\alpha Re)^{-1} i \int\limits_{K}\phi'_{n}\phi'_{k}dz \\
   & \quad
        + \int\limits_{K}U\phi_{n}\phi_{k}dz]
        + \sum_{K \in \ \mathcal{T}_{h}} \sum_{m \in \ I} p_{m_{\ast}}[
           \alpha^{-1} i \int\limits_{K}\psi_{m}\phi'_{k}dz] \\
   & \quad
        = c_{h} \cdot \{ \sum_{K \in \ \mathcal{T}_{h}} \sum_{n \in \ J} v_{n_{\ast}}
           \int\limits_{K}\phi_{n}\phi_{k}dz \}
\end{split}
\end{equation}

\begin{equation}
\label{eq_22}
\begin{split}
   & \sum_{K \in \ \mathcal{T}_{h}}
        \sum_{m \in \ I} p_{m_{\ast}}[
           - \alpha^{2} \int\limits_{K}\psi_{m}\psi_{\ell}dz
           - \int\limits_{K}\psi'_{m}\psi'_{\ell}dz] \\
   & \quad
        = \ - \sum_{K=K_{N}; \ n \in \ J_{K_{N}}; n \neq 3}
           v_{n_{\ast}} Re^{-1} \Phi_{n_{\ast}}''(a) \Psi_{N+1}(a)  \\
   & \quad
        + \sum_{K=K_{0}; \ n \in \ J_{K_{0}}; n \neq 1}
           v_{n_{\ast}} Re^{-1} \Phi_{n_{\ast}}''(0) \Psi_{1}(0)  \\
   & \quad
           + \sum_{K \in \ \mathcal{T}_{h}}
              \sum_{n \in \ J} v_{n_{\ast}}(
                 - 2 \alpha i \int\limits_{K}U'\phi_{n}\psi_{\ell}dz),
\end{split}
\end{equation}
for all $k \in J$, $\ell \in I$, for all $K \in \
\mathcal{T}_{h}$, $k \neq 1$ when $K=K_{0}$, $k \neq 3$ when
$K=K_{N}$.

The eigenvalue problem in $c_{h}$ $\in$ $\mathbb{C}$, $A_{h}=($
$u_{1},v_{1},$ $\ldots,$ $u_{2N+1},v_{2N+1})$ $\in$
$\mathbb{C}^{2(2N+1)}$, $B_{h}=($ $p_{0},$ $p_{1},$ $\ldots,$
$p_{N+1})$ $\in$ $\mathbb{C}^{N+2}$ represented by relations
$(\ref{eq_20})$ - $(\ref{eq_22})$ can be written in the following
matrix form
\begin{equation}
\label{NS5_p}
   K_{h}A_{h} + L_{h}B_{h} = c_{h}S_{h}A_{h},
      \quad G_{h}B_{h} = H_{h}A_{h} ,
\end{equation}
where $K_{h}, S_{h}$ $\in$ $\mathbb{C}^{2(2N+1) \times 2(2N+1)}$,
$L_{h}$ $\in$ $\mathbb{C}^{2(2N+1) \times (N+2)}$, $G_{h}$ $\in$
$\mathbb{C}^{(N+2) \times (N+2)}$, $H_{h}$ $\in$
$\mathbb{C}^{(N+2) \times 2(2N+1)}$.

Matrix $G_{h}$ is positive definite and symmetric. So problem
(\ref{NS5_p}) becomes
\begin{equation}
\label{NS5_pp}
   (K_{h} + L_{h}G_{h}^{-1}H_{h})A_{h} = c_{h}S_{h}A_{h},
      \quad B_{h} = G_{h}^{-1}H_{h}A_{h} ,
\end{equation}
The eigenvalue problem given by the first equation of
(\ref{NS5_pp}) can be solved by using the QZ method or, after
multiplication by $S_{h}^{-1}$ at the left, by using the QR or LR
method.

\section{Possible applications to concrete fluid flows}
\label{sectiunea4}

Once $(\ref{NS5_pp})$ solved, several linear hydrodynamic
stability characteristics can be determined numerically. Keeping
$\alpha$ and $Re$ fixed, (\ref{NS5_pp}) yields approximations of
$c$ and of amplitude distributions $u(y)$, $v(y)$ and $p(y)$.
Repeating this operation with various values of $\alpha$ and $Re$,
we can identify the pairs ($Re,\alpha$) where the approximate
value $c_{i}=0$. Thus, we can construct the neutral curve. We also
can determine the curves $c_{i}(Re,\alpha)=constant$ of constant
amplification factor and the wave speed $c_{r}(Re,\alpha)$, for
$(Re,\alpha)$ belonging to the neutral curve.

In particular, problem $(\ref{NS5_pp})$ is the basis of a computer
program for the case of Prandtl's boundary layer, which will be
presented elsewhere.

\end{document}